\documentclass[11pt,twoside]{article}

\usepackage[latin1]{inputenc}
\usepackage{epsfig}
\usepackage{amsfonts}
\usepackage{cite}
\usepackage{color}
\definecolor{lgrey}{gray}{0.99}

\title{
  {\huge Complex Analysis of Real Functions \\[1.5ex]}
  VII: A Simple Extension of the \\
  Cauchy-Goursat Theorem }

\author{
  \Large Jorge L. deLyra\footnote{Email: delyra@latt.if.usp.br} \\
  Department of Mathematical Physics \\
  Physics Institute \\
  University of São Paulo }

\date{August 19, 2018}

\newcommand{\ii}{\mbox{\boldmath$\imath$}}

\newcommand{\ldot}{\mbox{\Large$\cdot$}\!}

\newtheorem{definition}{Definition}
\newtheorem{property}{Property}[definition]
\newtheorem{lemma}{Lemma}
\newtheorem{theorem}{Theorem}

\newtheorem{proof}{Proof}[theorem]
\newcommand{\Colon}{{\hspace{-0.4em}\bf:}}

\begin{document}\maketitle

\vspace{-3.4ex}
\begin{abstract}
  \noindent
  In the context of the complex-analytic structure within the open unit
  disk, that was established in a previous paper, here we establish a
  simple generalization of the Cauchy-Goursat theorem of complex analytic
  functions. We do this first for the case of inner analytic functions,
  and then generalize the result to all analytic functions. We thus show
  that the Cauchy-Goursat theorem holds even if the complex function has
  isolated singularities located {\em on} the integration contour, so long
  as these are all {\em integrable} ones.
\end{abstract}

\section{Introduction}\label{Sec01}

In previous papers~\cite{CAoRFI,CAoRFII,CAoRFIII,CAoRFIV,CAoRFV,CAoRFVI}
we have shown that there is a correspondence between, on the one hand,
real functions and other real objects on the unit circle, and on the other
hand, inner analytic functions within the open unit disk of the complex
plane~\cite{CVchurchill}. This correspondence is based on the
complex-analytic structure which we introduced in~\cite{CAoRFI}. That
complex-analytic structure includes the concept of inner analytic
functions, two analytic operations on them, which we named angular
differentiation and angular integration, and a scheme for the
classification of all the possible singularities of these functions.

This classification scheme separates the singularities as either soft or
hard ones, depending on whether or not the limit of the function to the
singular point exists. As part of that classification scheme we also
introduced gradations of both hardness and softness for the singularities,
given by integers degrees. In particular, a hard singularity which becomes
soft under a single angular integration of the inner analytic function is
a borderline hard one, with degree of hardness zero. In a previous
paper~\cite{CAoRFI} we have shown that both soft and borderline hard
singularities are integrable ones, while the hard singularities with
strictly positive degrees of hardness are non-integrable ones.

Here we will show that one can extend the Cauchy-Goursat theorem, by
weakening its hypotheses, so that both soft and borderline hard isolated
singularities are allowed on the integration contour. We will prove this
first for inner analytic functions, and then generalize the result to
arbitrary complex analytic functions, using conformal transformations.

For ease of reference, we include here a one-page synopsis of the
complex-analytic structure introduced in~\cite{CAoRFI}. It consists of
certain elements within complex analysis~\cite{CVchurchill}, as well as of
their main properties.

\paragraph{Synopsis:} The Complex-Analytic Structure\\

\noindent
An {\em inner analytic function} $w(z)$ is simply a complex function which
is analytic within the open unit disk. An inner analytic function that has
the additional property that $w(0)=0$ is a {\em proper inner analytic
  function}. The {\em angular derivative} of an inner analytic function is
defined by

\noindent
\begin{equation}
  w^{\ldot}(z)
  =
  \ii
  z\,
  \frac{dw(z)}{dz}.
\end{equation}

\noindent
By construction we have that $w^{\ldot}(0)=0$, for all $w(z)$. The {\em
  angular primitive} of an inner analytic function is defined by

\begin{equation}\label{EQAngInt}
  w^{-1\ldot}(z)
  =
  -\ii
  \int_{0}^{z}dz'\,
  \frac{w(z')-w(0)}{z'}.
\end{equation}

\noindent
By construction we have that $w^{-1\ldot}(0)=0$, for all $w(z)$. In terms
of a system of polar coordinates $(\rho,\theta)$ on the complex plane,
these two analytic operations are equivalent to differentiation and
integration with respect to $\theta$, taken at constant $\rho$. These two
operations stay within the space of inner analytic functions, they also
stay within the space of proper inner analytic functions, and they are the
inverses of one another. Using these operations, and starting from any
proper inner analytic function $w^{0\ldot}(z)$, one constructs an infinite
{\em integral-differential chain} of proper inner analytic functions,

\begin{equation}
  \left\{
    \ldots,
    w^{-3\ldot}(z),
    w^{-2\ldot}(z),
    w^{-1\ldot}(z),
    w^{0\ldot}(z),
    w^{1\ldot}(z),
    w^{2\ldot}(z),
    w^{3\ldot}(z),
    \ldots\;
  \right\}.
\end{equation}

\noindent
Two different such integral-differential chains cannot ever intersect each
other. There is a {\em single} integral-differential chain of proper inner
analytic functions which is a constant chain, namely the null chain, in
which all members are the null function $w(z)\equiv 0$.

A general scheme for the classification of all possible singularities of
inner analytic functions is established. A singularity of an inner
analytic function $w(z)$ at a point $z_{1}$ on the unit circle is a {\em
  soft singularity} if the limit of $w(z)$ to that point exists and is
finite. Otherwise, it is a {\em hard singularity}. Angular integration
takes soft singularities to other soft singularities, and angular
differentiation takes hard singularities to other hard singularities.

Gradations of softness and hardness are then established. A hard
singularity that becomes a soft one by means of a single angular
integration is a {\em borderline hard} singularity, with degree of
hardness zero. The {\em degree of softness} of a soft singularity is the
number of angular differentiations that result in a borderline hard
singularity, and the {\em degree of hardness} of a hard singularity is the
number of angular integrations that result in a borderline hard
singularity. Singularities which are either soft or borderline hard are
integrable ones. Hard singularities which are not borderline hard are
non-integrable ones.

Given an integrable real function $f(\theta)$ on the unit circle, one can
construct from it a unique corresponding inner analytic function $w(z)$.
Real functions are obtained through the $\rho\to 1_{(-)}$ limit of the
real and imaginary parts of each such inner analytic function and, in
particular, the real function $f(\theta)$ is obtained from the real part
of $w(z)$ in this limit. The pair of real functions obtained from the real
and imaginary parts of one and the same inner analytic function are said
to be mutually Fourier-conjugate real functions.

Singularities of real functions can be classified in a way which is
analogous to the corresponding complex classification. Integrable real
functions are typically associated with inner analytic functions that have
singularities which are either soft or at most borderline hard. This ends
our synopsis.

\vspace{2.6ex}

\noindent
Unlike what was the case for the previous papers in this series, which
were focused mostly on the role played by this complex-analytic structure
in the analysis of real functions and other real objects, the result
presented in this paper concerns directly the theory of complex analytic
functions. Before we tackle the central issue of this paper, we must
present a review and refinement of some previous results, which were part
of that complex-analytic structure.

For the work to be developed in this paper it is important to recall that,
as was already noted in~\cite{CAoRFI}, whenever the singularities on the
unit circle are branch points, the corresponding branch cuts are to be
extended {\em outward\/} from the unit circle, so that there are no branch
cuts crossing the unit disk. In fact, in order to simplify the arguments
to be developed here, we adopt this recipe as an integral part of the
complex-analytic structure.

Some of the material contained in this paper can be seen as a development,
reorganization and extension of part of the material found, sometimes
still in rather rudimentary form, in the
papers~\cite{FTotCPI,FTotCPII,FTotCPIII,FTotCPIV,FTotCPV}.

\section{Refinement of Two Previous Results}\label{Sec02}

The discussions of two of the results related to the complex-analytic
structure introduced in~\cite{CAoRFI}, whose proofs were presented in that
paper, turn out to be somewhat incomplete for our needs here, so that the
proofs presented there must be somewhat refined. These are the discussions
regarding the fact that soft singularities must be integrable ones
(Property~$5.1$ of Definition~$5$), and the fact that borderline hard
singularities must be integrable ones (Property~$8.1$ of Definition~$8$).
Let us discuss and refine each one in turn.

With regard to the fact that a soft singularity of an inner analytic
function $w(z)$ at a point $z_{1}$ on the unit circle must be an
integrable one, the discussion given in~\cite{CAoRFI} takes us to the
point where it is shown that the integral of $w(z)$ exists on all simple
curves contained within the unit disk that connect to the point $z_{1}$.
However, we neglected to point out that certain integrals involving $w(z)$
have all the same {\em value}, which is equivalent to the fact that the
angular primitive $w^{-1\ldot}(z)$ is well-defined at $z_{1}$. Let us
repeat the argument here. Therefore, we now review the following important
property of soft singularities, first established in~\cite{CAoRFI}.

\setcounter{definition}{5}\setcounter{property}{0}

\begin{property}\Colon\label{Prop1SoHaSi}
  A soft singularity of an inner analytic function $w(z)$ at a point
  $z_{1}$ on the unit circle is necessarily an integrable one.
\end{property}

\setcounter{definition}{0}\setcounter{property}{0}

\begin{figure}[t]
  \centering
  {\color{white}\rule{\textwidth}{0.1ex}}
  \fbox{
    \epsfig{file=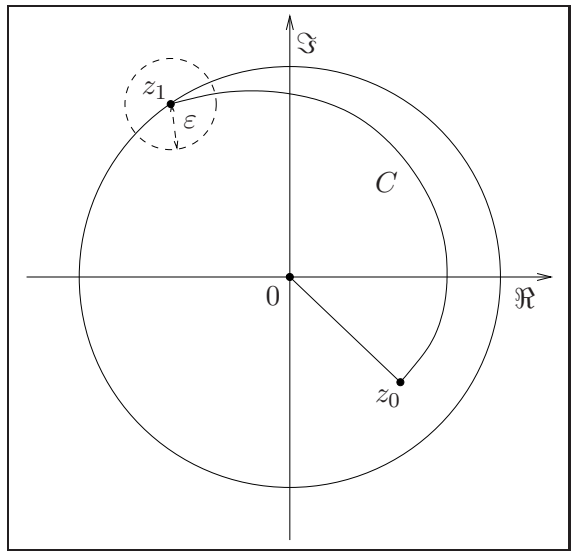,scale=1.0,angle=0}
  }
  \caption{The unit circle of the complex plane, the singular point
    $z_{1}$, the origin $0$, the internal point $z_{0}$, the neighborhood
    of radius $\varepsilon$, the integration contour $C$ and the straight
    segment connecting $z_{0}$ to the origin.}
  \label{Fig01}
\end{figure}

\noindent
This is so because the angular integration of $w(z)$ produces its angular
primitive, an inner analytic function $w^{-1\ldot}(z)$ which also has at
$z_{1}$ a soft singularity, and therefore is well defined at that point.
Since the value of $w^{-1\ldot}(z)$ at $z_{1}$ is given by an integral
involving $w(z)$ along a curve from the origin to $z_{1}$, as shown in
Equation~(\ref{EQAngInt}), that integral must therefore exist and result
in a finite complex number, for all curves within the open unit disk that
go from $0$ to $z_{1}$. Since the other factor involved in the integrand
of that integral is a regular function which is different from zero in a
neighborhood around $z_{1}$, say an open disk of radius $\varepsilon$ as
shown in Figure~\ref{Fig01}, this implies that $w(z)$ must be integrable
around $z_{1}$. Therefore, the singularity of $w(z)$ at $z_{1}$ must be an
integrable one.

In addition to this, since $w^{-1\ldot}(z)$ has a definite finite complex
value at $z_{1}$, we may also conclude that the value of the integral
giving it does not depend on the integration contour from $0$ to $z_{1}$,
that is, on the direction along which that curve connects to $z_{1}$. Then
the Cauchy-Goursat theorem in its usual form, which implies that integrals
from the origin to any point $z_{0}$ within the open unit disk are
independent of the integration contours connecting those two points and
contained within the open unit disk, allows us to generalize the result,
in a straightforward way, from curves starting at the point $0$ to curves
starting at any internal point $z_{0}$ within the open unit disk, by
simply connecting the origin and $z_{0}$ by means of any simple curve, for
example the straight segment shown in Figure~\ref{Fig01}. We therefore
conclude that, given an integration contour $C$ contained within the unit
disk and going from $z_{0}$ to $z_{1}$, the integral

\begin{equation}
  \int_{C}dz\,
  \frac{w(z)-w(0)}{z}
\end{equation}

\noindent
does not depend on the contour. We thus establish this property in a more
complete way.

\vspace{2.6ex}

\noindent
With regard to the fact that a borderline hard singularity of an inner
analytic function $w(z)$ at a point $z_{1}$ on the unit circle must be an
integrable one, we again neglected to point out in~\cite{CAoRFI} that the
same set of integrals from $z_{0}$ to $z_{1}$ discussed in the previous
case, along any simple curve contained within the unit disk and connecting
those two points, are all equal, which once more is equivalent to the fact
that the angular primitive $w^{-1\ldot}(z)$ is well-defined at $z_{1}$.
Let us repeat the argument here. Therefore, we now review the following
important property of borderline hard singularities, first established
in~\cite{CAoRFI}.

\setcounter{definition}{8}\setcounter{property}{0}

\begin{property}\Colon\label{Prop1GraSin}
  A borderline hard singularity of an inner analytic function $w(z)$ at a
  point $z_{1}$ on the unit circle must be an integrable one.
\end{property}

\setcounter{definition}{0}\setcounter{property}{0}

\noindent
This is so because the angular integration of $w(z)$ produces its angular
primitive, an inner analytic function $w^{-1\ldot}(z)$ which has at
$z_{1}$ a soft singularity, given that the singularity of $w(z)$ at that
point is a borderline hard one, and therefore $w^{-1\ldot}(z)$ is well
defined at the point $z_{1}$. Since the value of $w^{-1\ldot}(z)$ at
$z_{1}$ is given by an integral involving $w(z)$ along a curve from the
origin to $z_{1}$, as shown in Equation~(\ref{EQAngInt}), that integral
must therefore exist and result in a finite complex number, for all curves
within the unit disk that go from $0$ to $z_{1}$. Since the other factor
involved in the integrand of that integral is a regular function which is
different from zero in the neighborhood around $z_{1}$, this implies that
$w(z)$ must be integrable around $z_{1}$. Therefore, the singularity of
$w(z)$ at $z_{1}$ must be an integrable one.

In addition to this, since $w^{-1\ldot}(z)$ has a definite finite complex
value at $z_{1}$, we may also conclude that the value of the integral
giving it does not depend on the integration contour from the origin to
$z_{1}$, that is, on the direction along which that curve connects to
$z_{1}$. Just as was noted before in the discussion of the previous case,
at this point the Cauchy-Goursat theorem allows us to generalize the
result, in a straightforward way, from curves starting at the point $0$ to
curves starting at any internal point $z_{0}$ on the open unit disk. We
therefore conclude that, given an integration contour $C$ contained within
the unit disk and going from $z_{0}$ to $z_{1}$, the integral

\begin{equation}
  \int_{C}dz\,
  \frac{w(z)-w(0)}{z}
\end{equation}

\noindent
does not depend on the contour. We thus establish this property in a more
complete way.

\vspace{2.6ex}

\noindent
We may express these two results, about the singularities being
integrable, in a more concise way, by simply stating that what we mean by
the integrability of an inner analytic function $w(z)$ around a singular
point $z_{1}$ on the unit circle is that the integrals shown in
Equation~(\ref{EQAngInt}) on curves contained within the open unit disk
and going from any internal point $z_{0}$ to that singular point both
exist and are independent of the curves. This is, of course, equivalent to
the statement that the angular primitive of $w(z)$ exists and is
well-defined at $z_{1}$.

Once again we recall that, as was noted in~\cite{CAoRFI}, whenever the
singularities on the unit circle are branch points, it is understood that
the corresponding branch cuts are to be extended {\em outward\/} from the
unit circle, so that no branch cuts cross the unit disk. In this way there
can be no crossings between branch cuts and integration contours within
the open unit disk. This simplifies the arguments, since such crossings
would force us to consider the fact that the integration contours might be
changing from one leaf of a Riemann surface to another. However, the fact
that this is {\em not\/} an essential hypothesis is apparent when one
considers that for {\em closed\/} integration contours these crossings
would necessarily happen in pairs, each pair representing a change of
leafs followed by a change back to the original leaf, given that all
branch cuts within the open unit disk must cross it completely, since
there are no singularities of $w(z)$ within the open unit disk.

\section{Extension of the Cauchy-Goursat Theorem}\label{Sec03}

In either one of the two situations examined in Section~\ref{Sec02}, in
which the inner analytic function $w(z)$ has either a soft singularity or
a borderline hard singularity at a point $z_{1}$ on the unit circle, we
discovered that, given an integration contour $C$ contained within the
unit disk and going from any internal point $z_{0}$ to $z_{1}$, the
integral

\begin{equation}
  \int_{C}dz\,
  \frac{w(z)-w(0)}{z}
\end{equation}

\noindent
does not depend on the contour. We now observe that, if $w(z)$ is any
inner analytic function, then $w_{p}(z)=w(z)-w(0)$ is a {\em proper} inner
analytic function, since $w_{p}(0)=0$. Therefore what we have here is the
statement that the integral

\begin{equation}\label{EQPrInAnIndep}
  \int_{C}dz\,
  \frac{w_{p}(z)}{z}
\end{equation}

\noindent
does not depend on the contour, for all integration contours $C$ contained
within the unit disk that go from $z_{0}$ to $z_{1}$, and for all proper
inner analytic functions within the unit disk that have an integrable
singularity at $z_{1}$. Let us now consider the integral

\begin{equation}
  \int_{C}dz\,
  w(z),
\end{equation}

\noindent
for an arbitrary inner analytic function $w(z)$ that has an integrable
singularity at $z_{1}$. Since $w(z)$ is necessarily regular at $z=0$, it
follows that the function $w_{p}(z)=zw(z)$ is a {\em proper} inner
analytic function, given that $w_{p}(0)=0$. In addition to this, since the
function $z$ is analytic everywhere, it also follows that $w_{p}(z)$ and
$w(z)$ have the same singularity structure, and thus we can write

\begin{equation}
  \int_{C}dz\,
  w(z)
  =
  \int_{C}dz\,
  \frac{w_{p}(z)}{z},
\end{equation}

\noindent
where by the statement involving Equation~(\ref{EQPrInAnIndep}) this last
integral is independent of the integration contour $C$. Therefore, we have
the statement that

\begin{equation}
  \int_{C}dz\,
  w(z)
\end{equation}

\noindent
is independent of the contour, for all integration contours $C$ within the
unit disk that go from $z_{0}$ to $z_{1}$, and for all inner analytic
functions that have an integrable singularity at $z_{1}$.

\begin{figure}[t]
  \centering
  {\color{white}\rule{\textwidth}{0.1ex}}
  \fbox{
    \epsfig{file=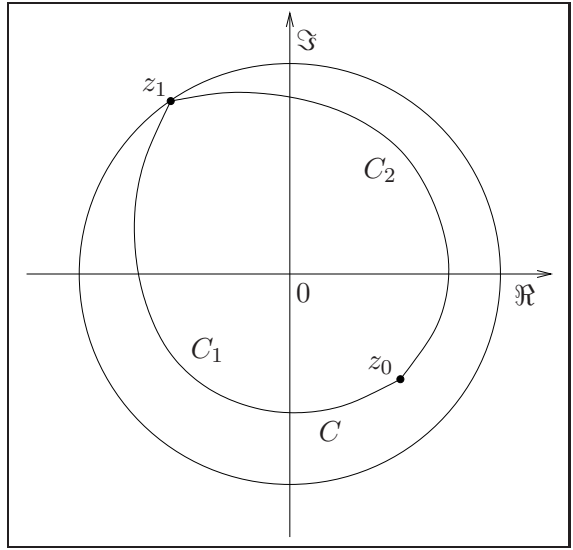,scale=1.0,angle=0}
  }
  \caption{The unit circle of the complex plane, the singular point
    $z_{1}$, the internal point $z_{0}$, and the integration contours
    involved in the proof of the extended Cauchy-Goursat theorem for inner
    analytic function within the unit disk.}
  \label{Fig02}
\end{figure}

We have thus determined that in these two cases the integral of an
arbitrary inner analytic function $w(z)$ from an arbitrary point $z_{0}$,
internal to the open unit disk, to the point $z_{1}$ on the unit circle,
where in either case $w(z)$ has an isolated integrable singularity, along
an integration contour contained within the closed unit disk and that
touches the unit circle only at $z_{1}$, is independent of that
integration contour from $z_{0}$ to $z_{1}$. Therefore, given two
different such curves, such as the curves $C_{1}$ and $C_{2}$ illustrated
in Figure~\ref{Fig02}, we may immediately conclude that the integral of
$w(z)$ over the {\em closed\/} integration contour $C$ formed by the two
curves is zero,

\begin{equation}\label{EQCauInnSing}
  \oint_{C}dz\,
  w(z)
  =
  0.
\end{equation}

\noindent
Since $z_{0}$ is an arbitrary internal point, this is valid for all closed
simple curves $C$ within the unit disk, that touch the unit circle only at
$z_{1}$. Observe that what we have concluded here is, in fact, that the
validity of the Cauchy-Goursat theorem, for the case of inner analytic
functions, is not disturbed by the presence of an isolated singularity
{\em on} the integration contour, so long as this singularity is either a
{\em soft\/} one or a {\em borderline hard\/} one or, in other words, so
long as singularities like this are all {\em integrable} ones. It is quite
clear, therefore, that this result constitutes an extension of the usual
form of the Cauchy-Goursat theorem, one which is valid at least for inner
analytic functions within the unit disk.

\begin{figure}[t]
  \centering
  {\color{white}\rule{\textwidth}{0.1ex}}
  \fbox{
    \epsfig{file=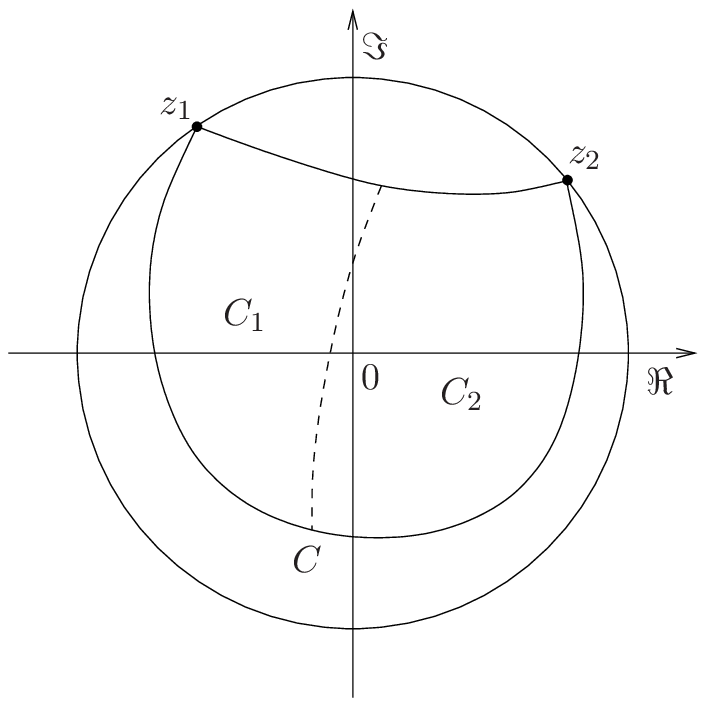,scale=1.0,angle=0}
  }
  \caption{The unit circle of the complex plane, the integration contour
    $C$ (solid line) with two singular points $z_{1}$ and $z_{2}$, showing
    how it can be decomposed into two integration contours $C_{1}$ and
    $C_{2}$, each one with only one singular point, by a simple cut
    (dashed line).}
  \label{Fig03}
\end{figure}

This result for a single point of singularity can then be trivially
extended, by contour manipulation and the repeated use of the
Cauchy-Goursat theorem in its usual form, to integration contours that
touch the unit circle on a finite number of points, at all of which $w(z)$
has isolated integrable singularities. Therefore, as a side effect of the
arguments presented in Section~\ref{Sec02}, in this section we will prove
the following theorem.

\begin{theorem}\Colon\label{Theo01}
  Given an inner analytic function $w(z)$, and a closed integration
  contour $C$ contained within the closed unit disk, that touches the unit
  circle only at a set of points satisfying one of two conditions, either
  points where $w(z)$ is analytic, or points where $w(z)$ has isolated
  integrable singularities, of which there must be a finite number, it
  follows that the integral of $w(z)$ over the contour $C$ is zero,

  \begin{equation}
    \oint_{C}dz\,
    w(z)
    =
    0.
  \end{equation}

\end{theorem}

\begin{proof}\Colon
\end{proof}

\noindent
We start with contours that touch the unit circle at a single integrable
singular point $z_{1}$, such as the one shown in Figure~\ref{Fig02}, and
by just using the results of Section~\ref{Sec02}, as expressed by
Equation~(\ref{EQCauInnSing}), to simply state that the integral of $w(z)$
over any such contour is zero. Next, given a closed contour $C$ that
touches the unit circle at two separate singular points $z_{1}$ and
$z_{2}$, such as the one shown in Figure~\ref{Fig03}, it can always be
separated into two closed contours $C_{1}$ and $C_{2}$, each one of which
touches the unit circle at only one singular point, as in the example
shown in Figure~\ref{Fig03}, by a simple cut (dashed line). When the two
separate closed contours $C_{1}$ and $C_{2}$ are joined together to form
the complete contour $C$, due to the orientation of the contours the
integrals over the cut, which is traversed twice, once in each direction,
cancel out. This can also be done for a contour that touches the unit
circle at any finite number of separate singular points.

\begin{figure}[t]
  \centering
  {\color{white}\rule{\textwidth}{0.1ex}}
  \fbox{
    \epsfig{file=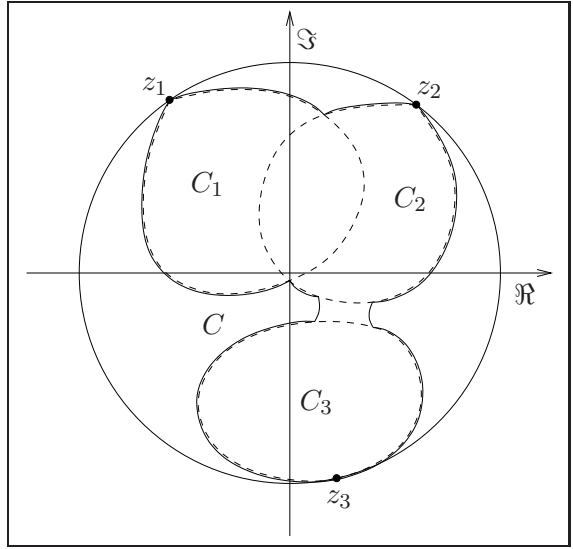,scale=1.0,angle=0}
  }
  \caption{The unit circle of the complex plane, the singular points
    $z_{1}$, $z_{2}$ and $z_{3}$, and the corresponding integration
    contours $C_{1}$, $C_{2}$ and $C_{3}$ (dashed lines), joined into a
    single overall contour $C$ (solid line).}
  \label{Fig04}
\end{figure}

Another way to think about this is to consider that one can construct any
contour such as those described in the statement of Theorem~\ref{Theo01}
by joining together a finite number of contours, each one of which touches
the unit circle at only one singular point, as is illustrated for the case
of three contours in Figure~\ref{Fig04}. Thus we see that the enormous
freedom to deform integration contours within the open unit disk without
changing the value of the integrals, which is given to us by the
Cauchy-Goursat theorem in its usual form, can be used to reduce a generic
closed contour, that touches $N$ isolated integrable singular points on
the unit circle, to a set of $N$ closed contours, each one of which
touches the unit circle at only one such singular point. This effectively
reduces the proof for the large and more complex contour to that for the
simple contour with only one singularity.

In addition to this, using the Cauchy-Goursat theorem in its usual form,
we may also deform any contour so that it morphs into one that touches the
unit circle at any points on that circle where $w(z)$ is analytic. In
other words, the integration contour may also run along any parts of the
unit circle on which $w(z)$ has no singularities at all. This completes,
therefore, the proof of Theorem~\ref{Theo01}.

\vspace{2.6ex}

\noindent
In this section we have established that the extended version of the
Cauchy-Goursat theorem, allowing for the presence of a finite number of
isolated integrable singularities on the integration contour, holds for
all inner analytic functions within the unit disk. In Section~\ref{Sec05}
we will generalize that result to all complex analytic functions, anywhere
on the complex plane, using conformal transformations. Therefore, before
we discuss the generalization of the theorem we must discuss these
conformal transformations.

\section{Conformal Transformations and Singularities}\label{Sec04}

The validity of the extended version of the Cauchy-Goursat theorem can be
generalized to all complex analytic functions integrated on arbitrarily
given closed integration contours, through the use of conformal
transformations. In order to do this, let us first establish the
definition and the notation for a conformal transformation. This is
essentially a shorter version of the discussion on this topic which was
given in Section~4 of a previous paper~\cite{CAoRFV}. Consider therefore
two complex variables $z_{a}$ and $z_{b}$ and the corresponding complex
planes, a complex analytic function $\gamma(z)$ defined on the complex
plane $z_{a}$ with values on the complex plane $z_{b}$, and its inverse
function, which is a complex analytic function $\gamma^{(-1)}(z)$ defined
on the complex plane $z_{b}$ with values on the complex plane $z_{a}$,

\noindent
\begin{eqnarray}
  z_{b}
  & = &
  \gamma(z_{a}),
  \nonumber\\
  z_{a}
  & = &
  \gamma^{(-1)}(z_{b}).
\end{eqnarray}

\noindent
Let us point out here that these relations immediately imply that

\noindent
\begin{eqnarray}\label{EQTransfDerivs}
  \frac{dz_{b}}{dz_{a}}
  & = &
  \frac{d\gamma(z_{a})}{dz_{a}},
  \nonumber\\
  \frac{dz_{a}}{dz_{b}}
  & = &
  \frac{d\gamma^{(-1)}(z_{b})}{dz_{b}},
\end{eqnarray}

\noindent
which in turn imply that

\begin{equation}
  \frac{d\gamma(z_{a})}{dz_{a}}\,
  \frac{d\gamma^{(-1)}(z_{b})}{dz_{b}}
  =
  1,
\end{equation}

\noindent
for all pairs of points $z_{a}$ and $z_{b}$ related by the conformal
transformation. This means that any point where the derivative of
$\gamma(z_{a})$ has a zero on the $z_{a}$ plane corresponds to a point
where the derivative of $\gamma^{(-1)}(z_{b})$ has a singularity on the
$z_{b}$ plane, and vice-versa.

Consider a bounded and simply connected open region $S_{a}$ on the complex
plane $z_{a}$ and its image $S_{b}$ under $\gamma(z)$, which is a similar
region on the complex plane $z_{b}$. It can be shown that if $\gamma(z)$
is analytic on $S_{a}$, is invertible there, and its derivative has no
zeros there, then its inverse function $\gamma^{(-1)}(z)$ has these same
three properties on $S_{b}$, and the mapping between the two complex
planes established by $\gamma(z)$ and $\gamma^{(-1)}(z)$ is conformal, in
the sense that it preserves the angles between oriented curves at points
where they cross each other. The famous {\em Riemann mapping theorem}
states that such a conformal transformation $\gamma(z)$ always exists
between the open unit disk $S_{a}$ and any region $S_{b}$. In addition to
this, these properties of $\gamma(z)$ and $\gamma^{(-1)}(z)$ can be
extended to the boundary of the regions so long as these boundaries are
differentiable simple curves.

Consider therefore that the regions under consideration are the interiors
of simple closed curves. One of these curves will be the unit circle
$C_{a}$ on the complex plane $z_{a}$, and the other will be a given closed
differentiable simple curve $C_{b}$ on the complex plane $z_{b}$. Since
$\gamma(z_{a})$, being analytic, is in particular a continuous function,
the image on the $z_{b}$ plane of the unit circle $C_{a}$ on the $z_{a}$
plane must be a continuous closed curve $C_{b}$. We can also see that
$C_{b}$ must be a simple curve, because the fact that $\gamma(z_{a})$ is
invertible on $C_{a}$ means that it cannot have the same value at two
different points of $C_{a}$, and therefore no two points of $C_{b}$ can be
the same. Consequently, the curve $C_{b}$ cannot self-intersect. Finally,
the fact that $C_{b}$ must be a differentiable curve is a simple
consequence of the facts that the $\gamma(z_{a})$ transformation is
conformal and that the unit circle $C_{a}$ is a differentiable curve.

Given any analytic function $w_{a}(z_{a})$ on the $z_{a}$ plane, the
conformal transformation $\gamma(z_{a})$ maps it to a corresponding
function $w_{b}(z_{b})$ on the $z_{b}$ plane, and the inverse conformal
transformation $\gamma^{(-1)}(z_{b})$ maps that function back to the
function $w_{a}(z_{a})$ on the $z_{a}$ plane. We can do this by simply
composing either $w_{b}(z_{b})$ or $w_{a}(z_{a})$ with either the
transformation or its inverse, and simply passing the values of the
functions,

\noindent
\begin{eqnarray}\label{EQTransfFuncts}
  w_{b}(z_{b})
  & = &
  w_{a}(z_{a})
  \nonumber\\
  & = &
  w_{a}\!\left(\gamma^{(-1)}(z_{b})\right),
  \nonumber\\
  w_{a}(z_{a})
  & = &
  w_{b}(z_{b})
  \nonumber\\
  & = &
  w_{b}\!\left(\gamma(z_{a})\right).
\end{eqnarray}

\noindent
Since the composition of two analytic functions is also an analytic
function, and since $\gamma(z_{a})$ is analytic on the closed unit disk,
whenever $w_{b}(z_{b})$ is analytic on the $z_{b}$ plane the corresponding
function $w_{a}(z_{a})$ will also be analytic at the corresponding points
on the $z_{a}$ plane. Of course, where one of these two functions has an
isolated singularity on its plane of definition, so will the other on the
corresponding point in the other plane. Note that if any of these
functions is integrated over a closed integration contour on which it has
any isolated integrable singularities, whenever these singularities are
branch points we assume that the corresponding branch cuts extend {\em
  outward\/} from the integration contours.

The concepts of a soft singularity and of a borderline hard singularity
can be immediately extended from the case of inner analytic functions
within the unit disk to singularities of arbitrary complex analytic
functions anywhere on the complex plane. The concept of a soft singularity
of $w(z)$ at $z_{1}$ depends only on the existence of the $z\to z_{1}$
limit of $w(z)$. The concept of a hard singularity of $w(z)$ at $z_{1}$
depends only on the non-existence of that same limit. Finally, the concept
of a borderline hard singularity can be defined as that of a hard
singularity which is nevertheless an integrable one. In order to discuss
what happens with the singularities under the conformal transformation, we
must establish a few simple preliminary results, by means of the following
lemmas.

\begin{lemma}\Colon\label{Lemma01}
  If $z_{b,1}$ is a singular point of $w_{b}(z_{b})$, then the
  corresponding point $z_{a,1}$ under the conformal transformation is a
  singular point of $w_{a}(z_{a})$.
\end{lemma}

\noindent
Since according to Equations~(\ref{EQTransfFuncts}) we have that
$w_{b}(z_{b,1})=w_{a}\!\left(\hspace{-1.4pt}\gamma^{(-1)}(z_{b,1})\right)$
and since $\gamma^{(-1)}(z_{b})$ is analytic within and on the image of
the unit circle by the conformal transformation, if $w_{a}(z_{a})$ were
analytic at $z_{a,1}$, then $w_{b}(z_{b})$ would be analytic at $z_{b,1}$,
because the composition of two analytic functions is also an analytic
function. Therefore, if $w_{b}(z_{b})$ is singular at $z_{b,1}$, then
$w_{a}(z_{a})$ must be singular at $z_{a,1}$. This establishes
Lemma~\ref{Lemma01}.

\vspace{2.6ex}

\begin{lemma}\Colon\label{Lemma02}
  If the singularity of $w_{b}(z_{b})$ at $z_{b,1}$ is a soft one, then
  the singularity of $w_{a}(z_{a})$ at the corresponding singular point
  $z_{a,1}$ under the conformal transformation is also a soft singularity.
\end{lemma}

\noindent
Since according to the definition in Equations~(\ref{EQTransfFuncts}) we
have that $w_{a}(z_{a})=w_{b}(z_{b})$ and since the $z_{b}\to z_{b,1}$
limit on the $z_{b}$ plane corresponds, through the continuous conformal
transformation, to the $z_{a}\to z_{a,1}$ limit on the $z_{a}$ plane, it
follows that if the limit

\begin{equation}
  \lim_{z_{b}\to z_{b,1}}
  w_{b}(z_{b})
\end{equation}

\noindent
exists, then so does the limit

\begin{equation}
  \lim_{z_{a}\to z_{a,1}}
  w_{a}(z_{a}).
\end{equation}

\noindent
Therefore, if the singularity of $w_{b}(z_{b})$ at $z_{b,1}$ is a soft
one, which means that the first limit exists, then the singularity of
$w_{a}(z_{a})$ at $z_{a,1}$ is also a soft singularity, since in this case
the second limit also exists. This establishes Lemma~\ref{Lemma02}.

\vspace{2.6ex}

\begin{lemma}\Colon\label{Lemma03}
  If the singularity of $w_{b}(z_{b})$ at $z_{b,1}$ is a hard one, then
  the singularity of $w_{a}(z_{a})$ at the corresponding singular point
  $z_{a,1}$ under the conformal transformation is also a hard singularity.
\end{lemma}

\noindent
By an argument similar to the one used in Lemma~\ref{Lemma02}, since
according to the definition in Equations~(\ref{EQTransfFuncts}) we have
that $w_{a}(z_{a})=w_{b}(z_{b})$, and since the $z_{b}\to z_{b,1}$ and
$z_{a}\to z_{a,1}$ limits correspond to one another, if $w_{b}(z_{b})$ is
{\em not\/} well defined at $z_{b,1}$, which means that the limit

\begin{equation}
  \lim_{z_{b}\to z_{b,1}}
  w_{b}(z_{b})
\end{equation}

\noindent
does {\em not\/} exist, then the limit

\begin{equation}
  \lim_{z_{a}\to z_{a,1}}
  w_{a}(z_{a})
\end{equation}

\noindent
also does {\em not\/} exist, and therefore $w_{a}(z_{a})$ {\em cannot\/}
be well defined at the corresponding point $z_{a,1}$. Therefore, if the
singularity of $w_{b}(z_{b})$ at $z_{b,1}$ is a hard one, then the
singularity of $w_{a}(z_{a})$ at $z_{a,1}$ must also be a hard
singularity. This establishes Lemma~\ref{Lemma03}.

\vspace{2.6ex}

\begin{lemma}\Colon\label{Lemma04}
  If the singularity of $w_{b}(z_{b})$ at $z_{b,1}$ is an integrable one,
  then the singularity of $w_{a}(z_{a})$ at the corresponding singular
  point $z_{a,1}$ under the conformal transformation is also an integrable
  singularity.
\end{lemma}

\noindent
If the singularity of $w_{b}(z_{b})$ at $z_{b,1}$ is an integrable one,
then the integral of $w_{b}(z_{b})$ along any open integration contour
$D_{b}$ going from a point $z_{b,0}$ internal to $C_{b}$ to the singular
point $z_{b,1}$,

\begin{equation}
  \int_{D_{b}}dz_{b}\,
  w_{b}(z_{b}),
\end{equation}

\noindent
exists and is a finite complex number. We consider now the corresponding
integral of $w_{a}(z_{a})$ along an arbitrary open integration contour
$D_{a}$ going from an internal point $z_{a,0}$ within the open unit disk
to the singular point $z_{a,1}$ on the unit circle, and make a
transformation of the integration variable from $z_{a}$ to $z_{b}$,

\noindent
\begin{eqnarray}
  \int_{D_{a}}dz_{a}\,
  w_{a}(z_{a})
  & = &
  \int_{D_{b}}dz_{b}\,
  \left(
    \frac{dz_{a}}{dz_{b}}
  \right)
  w_{b}(z_{b})
  \nonumber\\
  & = &
  \int_{D_{b}}dz_{b}\,
  \left[
    \frac{d\gamma^{(-1)}(z_{b})}{dz_{b}}
  \right]
  w_{b}(z_{b}),
\end{eqnarray}

\noindent
where we used the relations shown in Equations~(\ref{EQTransfDerivs})
and~(\ref{EQTransfFuncts}), where $D_{b}$ is an open integration contour
going from an internal point $z_{b,0}$ to the singular point $z_{b,1}$,
and where $D_{b}$, $z_{b,0}$ and $z_{b,1}$ correspond respectively to
$D_{a}$, $z_{a,0}$ and $z_{a,1}$, through the conformal transformation.
Since the conformal transformation is an analytic function, the derivative
within brackets is also an analytic function within and on $D_{b}$, and
therefore is a limited function there. Since $w_{b}(z_{b})$ by hypothesis
is an integrable function around the singular point $z_{b,1}$, and since
the product of a limited function and an integrable function is also an
integrable function, the integrand of this last integral is an integrable
function, and therefore this last integral exists and is a finite complex
number. We thus conclude that

\begin{equation}
  \int_{D_{a}}dz_{a}\,
  w_{a}(z_{a})
\end{equation}

\noindent
exists and is a finite complex number. Therefore, $w_{a}(z_{a})$ is an
integrable function around the singular point $z_{a,1}$, and therefore
that singularity is also an integrable one. This establishes
Lemma~\ref{Lemma04}.

\vspace{2.6ex}

\noindent
We must now consider the question of what is the set of curves $C_{b}$ for
which the structure described above can be set up. Given the curve
$C_{b}$, the only additional objects we need in order to do this is the
conformal mapping $\gamma(z_{a})$ and its inverse $\gamma^{(-1)}(z_{b})$,
between that curve and the unit circle $C_{a}$, as well as between the
respective interiors. The existence of these transformation functions can
be ensured as a consequence of the Riemann mapping theorem, and of the
associated results relating to conformal mappings between regions of the
complex plane~\cite{RMPQiu}. According to that theorem, a conformal
transformation such as the one we just described exists between any
bounded simply connected open set of the plane and the open unit disk, and
can be extended to the respective boundaries so long as the curve $C_{b}$
is differentiable. There are therefore no additional limitations on the
differentiable simple closed curves $C_{b}$ that may be considered here.

\section{Generalization to Arbitrary Analytic Functions}\label{Sec05}

The basic idea of the proof of the generalization of the extended
Cauchy-Goursat theorem will be to embed an arbitrarily given integration
contour on the $z_{b}$ plane, which must be a simple closed curve but may
or may not be a differentiable curve, into a region bounded by a {\em
  differentiable} simple closed curve, which is then mapped to the unit
circle on the $z_{a}$ plane by a conformal transformation. The embedding
will be such that all the isolated integrable singularities on the
original integration contour are mapped onto the unit circle. This will
then allow us to use the extended Cauchy-Goursat theorem for inner
analytic functions within the unit disk, which was established in
Section~\ref{Sec03}, to prove the generalized version of that theorem.
Therefore, in this section we will prove the following theorem.

\begin{theorem}\Colon\label{Theo02}
  Given an analytic function $w(z)$, and a closed integration contour $C$
  within which $w(z)$ is analytic, and on which $w(z)$ is analytic except
  for a finite number of isolated integrable singularities, it follows
  that the integral of $w(z)$ over the contour $C$ is zero,

  \begin{equation}
    \oint_{C}dz\,
    w(z)
    =
    0.
  \end{equation}

\end{theorem}

\noindent
This is the most complete generalization of the extended Cauchy-Goursat
theorem that we will consider here, the only relevant limitation being
that the number of isolated integrable singularities be finite.

\begin{proof}\Colon
\end{proof}

\noindent
Let $w_{b}(z_{b})$ be an analytic function within a closed simple curve
$D_{b}$ on the $z_{b}$ plane, and let the function $w_{b}(z_{b})$ also be
analytic almost everywhere {\em on} $D_{b}$, with the exception of a
finite number of isolated integrable singular points. It follows that, due
to Lemmas~\ref{Lemma01}--\ref{Lemma04}, the corresponding function
$w_{a}(z_{a})$ on the $z_{a}$ plane will have the same number of
singularities on it, which will also be isolated integrable singular
points. The curve $D_{b}$ may not be differentiable at some points,
including at some of the singular points. We will consider the integral of
$w_{b}(z_{b})$ over the integration contour $D_{b}$ on the $z_{b}$ plane,
which will then, of course, correspond to the integral of $w_{a}(z_{a})$
over a corresponding integration contour $D_{a}$ on the $z_{a}$ plane,
under the conformal transformation.

\begin{figure}[t]
  \centering
  {\color{white}\rule{\textwidth}{0.1ex}}
  \fbox{
    \epsfig{file=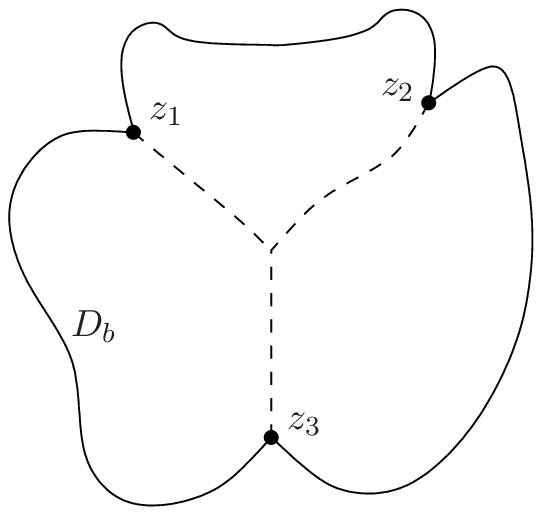,scale=1.0,angle=0}
  }
  \caption{An integration contour $D_{b}$ with three isolated
    singularities located at the points $z_{1}$, $z_{2}$ and $z_{3}$,
    where that contour is non-differentiable and concave, showing how to
    decompose it into three contours on which the singularities are
    located at points where the new contours are non-differentiable but
    convex, by the cuts shown (dashed lines).}
  \label{Fig05}
\end{figure}

The proof that follows will depend on the integration contour $D_{b}$
being either differentiable or non-differentiable and {\em convex} at all
the singular points found on it. However, this is not a true limitation,
because an integration contour that has one of more singular points where
it is non-differentiable and concave can always be decomposed into two or
more integration contours where those same singular points are convex, as
is shown in Figure~\ref{Fig05} for an example with three such points. As
one can see, all the three integration contours into which the original
one is decomposed by the cuts shown (dashed lines) are convex at the
singular points. When the three are put together to form once again the
original contour, the integrals over those cuts, which due to their
orientation are traversed once in each direction, cancel out. Therefore,
if the theorem is proven for all contours which are convex at the singular
points, it follows that it in fact holds for all contours, regardless of
whether they are convex or concave at their singular points. We may
therefore limit the proof to contours which are convex at all their
singular points.

\begin{figure}[t]
  \centering
  {\color{white}\rule{\textwidth}{0.1ex}}
  \fbox{
    \epsfig{file=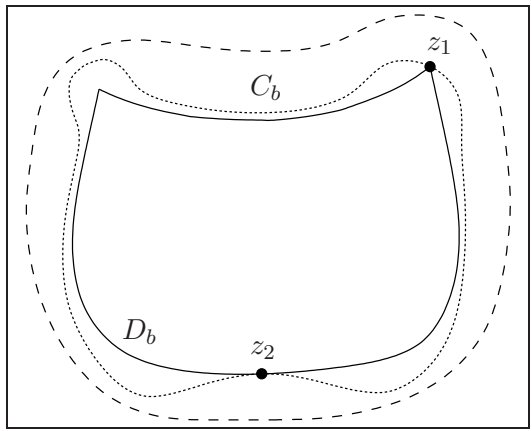,scale=1.0,angle=0}
  }
  \caption{The construction a differentiable closed curve $C_{b}$
    containing the integration contour $D_{b}$, to be conformally mapped
    to the unit circle $C_{a}$, showing also the singular points $z_{1}$
    and $z_{2}$, as well as a point of non-differentiability of the
    contour at which $w(z)$ is analytic.}
  \label{Fig06}
\end{figure}

Given an integration contour $D_{b}$ within which $w_{b}(z_{b})$ is
analytic, and on which $w_{b}(z_{b})$ is analytic except for a finite
number of isolated integrable singularities, at all of which the contour
is either differentiable or non-differentiable and convex, we consider now
the construction on the $z_{b}$ plane of a new closed differentiable
simple curve $C_{b}$ that contains $D_{b}$. At any points on $D_{b}$ where
$w_{b}(z_{b})$ is analytic there is a neighborhood of that point within
which there are no singularities of $w_{b}(z_{b})$ which are not located
directly on $D_{b}$, whose union forms a strip around $D_{b}$. In this
case we make $C_{b}$ go through these neighborhoods in a differentiable
fashion, outside the interior of $D_{b}$, so that we ensure that no
singularities of $w_{b}(z_{b})$ get included on $C_{b}$ or in its
interior, other than those on $D_{b}$. This can be done even at
non-singular points where the integration contour $D_{b}$ is {\em not\/}
differentiable, in which case we make $C_{b}$ just go around the point of
non-differentiability of $D_{b}$, in a differentiable fashion, as can be
seen illustrated in Figure~\ref{Fig06}.

At points of $D_{b}$ where $w_{b}(z_{b})$ has an isolated integrable
singularity, since the singularity is isolated, there is also a
neighborhood of the point within which there are no {\em other}
singularities of $w_{b}(z_{b})$, which is part of the aforementioned strip
around $D_{b}$. In this case we make $C_{b}$ go through this neighborhood,
still keeping to the outer side of $D_{b}$, in such a way that the curve
{\em runs over\/} the singular point in a differentiable way, which is
possible because $D_{b}$ is convex at that singular point, as is
illustrated by the point $z_{1}$ in Figure~\ref{Fig06}. The singular point
is one which the curve $C_{b}$ will therefore share with $D_{b}$, as is
also illustrated by the point $z_{1}$ in Figure~\ref{Fig06}. At singular
points where $D_{b}$ is differentiable we simply make $C_{b}$ tangent to
$D_{b}$ at that point, as is illustrated by the point $z_{2}$ in
Figure~\ref{Fig06}.

The result of this process, taken all around $D_{b}$ and including all the
isolated integrable singularities found on it, is a differentiable simple
curve $C_{b}$ which contains $D_{b}$ and the singularities on it, but that
contains no other singularities of $w_{b}(z_{b})$, and which shares with
$D_{b}$ all the points where the relevant isolated integrable
singularities of $w_{b}(z_{b})$ are located. Since by construction $C_{b}$
is a differentiable closed simple curve, by the Riemann mapping theorem
there exists a conformal transformation $\gamma(z_{a})$ that maps it from
the unit circle. Therefore, the inverse conformal transformation
$\gamma^{(-1)}(z_{b})$ will map all the isolated integrable singular
points on $D_{b}$ to the unit circle $C_{a}$.

Since the interior of $C_{b}$ in mapped by the inverse transformation
$\gamma^{(-1)}(z_{b})$ onto the open unit disk, it follows that the
integration contour $D_{b}$, which is contained in $C_{b}$, is mapped by
$\gamma^{(-1)}(z_{b})$ onto a closed simple integration contour $D_{a}$
contained in the unit disk in the $z_{a}$ plane, which will not be
differentiable if $D_{b}$ is not, but which is contained within the closed
unit disk, and that touches the unit circle only at each one of the
isolated integrable singular points of $w_{a}(z_{a})$ on $D_{a}$ that
correspond to the singularities of $w_{b}(z_{b})$ on $D_{b}$.

Observe that, since the curve $C_{b}$ does not contain any singularities
of $w_{b}(z_{b})$ in its strict interior, the interior of the curve
$C_{a}$, which is the open unit disk, does not contain any singularities
of $w_{a}(z_{a})$. Therefore, according to the definition given
in~\cite{CAoRFI}, $w_{a}(z_{a})$ is an inner analytic function. If we now
consider the integral of $w_{b}(z_{b})$ over $D_{b}$, it is a very simple
thing to change the integration variable from $z_{b}$ to $z_{a}$,

\noindent
\begin{eqnarray}
  \oint_{D_{b}}dz_{b}\,
  w_{b}(z_{b})
  & = &
  \oint_{D_{a}}dz_{a}\,
  \left(
    \frac{dz_{b}}{dz_{a}}
  \right)
  w_{a}(z_{a})
  \nonumber\\
  & = &
  \oint_{D_{a}}dz_{a}\,
  \left[
    \frac{d\gamma(z_{a})}{dz_{a}}
  \right]
  w_{a}(z_{a}),
\end{eqnarray}

\noindent
where we used the relations shown in Equations~(\ref{EQTransfDerivs})
and~(\ref{EQTransfFuncts}). Because $\gamma(z_{a})$ is an analytic
function on the whole closed unit disk, the derivative in brackets is also
an analytic function on the whole closed unit disk, and in addition to
this the function $w_{a}(z_{a})$ is analytic within the integration
contour $D_{a}$ and also on $D_{a}$ except for a finite set of isolated
singularities located on the unit circle. By the results of
Lemmas~\ref{Lemma01}--\ref{Lemma04}, these isolated singularity are all
integrable ones. Therefore, since the product of two analytic functions is
also an analytic function, the integrand is analytic within the
integration contour $D_{a}$, and also on it except for a finite set of
isolated integrable singularities on the unit circle, and hence is an
inner analytic function. Therefore, by Theorem~\ref{Theo01}, that is, the
extended Cauchy-Goursat theorem for inner analytic functions on the unit
disk, the integral is zero, and hence it follows that

\begin{equation}
  \oint_{D_{b}}dz_{b}\,
  w_{b}(z_{b})
  =
  0.
\end{equation}

\noindent
In other words, due to the fact that the integral of $w_{a}(z_{a})$ on
$D_{a}$ is zero, which is guaranteed by the extended Cauchy-Goursat
theorem for inner analytic functions, we may conclude that the integral of
$w_{b}(z_{b})$ on $D_{b}$ is also zero. This implies that the extended
Cauchy-Goursat theorem is valid for $w_{b}(z_{b})$, that is, for arbitrary
complex analytic functions anywhere on the complex plane. This completes
the proof of Theorem~\ref{Theo02}.

\vspace{2.6ex}

\noindent
Note that, once we have Theorem~\ref{Theo02} established, it is also valid
for all inner analytic functions, and therefore automatically includes the
contents of Theorem~\ref{Theo01}, which we may therefore regard as just an
intermediate step in the proof.

\section{Conclusions and Outlook}\label{Sec06}

We have shown that an extended version of the Cauchy-Goursat theorem holds
for all complex analytic functions, anywhere on the complex plane. The
extension of the theorem establishes that the integral of any such
function is zero, on any closed integration contour within which it is
analytic, even it the function has a finite set of isolated singularities
on the contour itself, so long as these are all integrable singularities.
It is interesting to note that this integrability requirement on the
singularities is the minimum necessary requirement for the integral over
the contour to make any sense at all, and it is very curious that it turns
out to be sufficient for the extension of the Cauchy-Goursat theorem.

The generalization of the result to infinite integration contours with a
countable infinity of isolated integrable singularities on them is quite
immediate, by a straightforward process of finite induction, using contour
manipulation, the extended Cauchy-Goursat theorem and the Cauchy-Goursat
theorem in its usual form. Therefore this extended Cauchy-Goursat theorem
can be used in essentially all circumstances in which the original
Cauchy-Goursat theorem applies, thus giving rise to extensions of many
previously known results.

Note that the establishment of this extended Cauchy-Goursat theorem gives
rise immediately to a corresponding extended set of Cauchy integral
formulas, when the interior of the integration contours contain isolated
poles. In fact, these extended Cauchy integral formulas are at the root of
the representation of integrable real functions by inner analytic
functions within the unit disk, which was established in~\cite{CAoRFI}. In
that paper the validity of these extended Cauchy integral formulas is
reflected by the fact that integrals which give the coefficients of the
Taylor series of the inner analytic functions are not only independent of
the radius $0<\rho<1$ of the circle over which they are calculated, but
are also continuous from within when one approaches the unit circle from
within the open unit disk, that is, in the $\rho\to 1_{(-)}$ limit.

Since the Cauchy-Goursat theorem is such an important and fundamental one,
it is to be expected that this extension will have further interesting
consequences, possibly in many fields of mathematics and also in many
applications in physics.

\section*{Acknowledgments}

The author would like to thank his friend and colleague Prof. Carlos
Eugênio Imbassay Carneiro, to whom he is deeply indebted for all his
interest and help, as well as his careful reading of the manuscript and
helpful criticism regarding this work.

\bibliography{allrefs_en}\bibliographystyle{ieeetr}

\end{document}